\documentclass [a4paper,10pt]{article}
\usepackage{amsmath}
\usepackage{amsfonts}
\usepackage{amssymb}
\usepackage{array,longtable}
\usepackage{amscd}
\usepackage[cp1251]{inputenc}
\usepackage[english]{babel}
\usepackage[usenames]{color}
\usepackage[matrix,arrow,curve]{xy}

\begin{document}
\newtheorem{definition}{Definition}
\newtheorem{theorem}{Theorem}
\newtheorem{lemma}{Lemma}
\newtheorem{proposition}{Proposition}
\newtheorem{conjecture}{Conjecture}
\newtheorem{criterion}{Criterion}
\newtheorem{observation}{Observation}
\newtheorem{corollary}{Corollary}
\newcommand{\diag}{{\rm diag}}
\newcommand{\tr}{{\rm Tr}}
\begin{center}
{\large {\bf HOW TO CHECK $D$-STABILITY: \linebreak A SIMPLE DETERMINANTAL TEST}} \\[0.5cm]
{\bf Olga Y. Kushel} \\[0.5cm]
Shanghai University, \\ Department of Mathematics, \\ Shangda Road 99, \\ 200444 Shanghai, China \\
kushel@mail.ru
\end{center}
\begin{abstract}

The concept of matrix $D$-stability, introduced in 1958 by Arrow and McManus is of major importance due to the variety of its applications. However, characterization of matrix $D$-stability for dimensions $n > 4$ is considered as a hard open problem. In this paper, we propose a simple way for testing matrix $D$-stability, in terms of the inequalities between principal minors of a matrix. The conditions are just sufficient but they allow to test matrices of an arbitrary size $n$, are easy to verify and can be used for the analysis of parameter-dependent models.

Hurwitz stability, $D$-stability, principal minors, $P$-matrices, determinantal inequalities.

MSC[2010] 15A12, 15A18, 15A75
\end{abstract}
 \section{Introduction}
 Let $\mathcal{M}^{n\times n}$ be
the set of all square real\ $n\times n$ matrices; $\sigma (\mathbf{A})$ be the
spectrum of a matrix $\mathbf{A}$ $\in \mathcal{M}^{n\times n}$ (i.e. the set of all eigenvalues of $\mathbf A$ defined as zeroes of its characteristic polynomial $f_{\mathbf A}(\lambda):= \det({\mathbf A} - \lambda{\mathbf I})$).

Let us recall the following crucial definitions.

{\bf Definition 1} (see, for example, \cite{BELL}, \cite{KAB}). A matrix $\mathbf{A}$ $\in \mathcal{M}^{n\times n}$ is
called {\it Hurwitz stable} or just {\it stable} if ${\rm Re}(\lambda )<0$ for all $\lambda \in \sigma (\mathbf{A})$.

{\bf Definition 2} (see \cite{AM}). A matrix $\mathbf{A}$ $\in \mathcal{M}^{n\times n}$ is
called {\it (multiplicative) $D$-stable} if ${\rm Re}(\lambda )<0$ for all $\lambda \in
\sigma (\mathbf{DA})$, where ${\mathbf D}$ is any positive diagonal matrix.

{\bf Definition 3} (see \cite{CROSS}). A matrix ${\mathbf A} \in {\mathcal M}^{n \times n}$ is called {\it diagonally stable} if the matrix $${\mathbf W}:={\mathbf D}{\mathbf A} + {\mathbf A}^{T}{\mathbf D}$$ is negative definite for some positive diagonal matrix ${\mathbf D}$.

Here, our purpose is to present new verifiable sufficient conditions for $D$-stability, that can be applied to an arbitrary $n \times n$ matrix. Though the problem of characterizing $D$-stability has been studied for over sixty years (see, e.g. \cite{KU2}, \cite{KAB} and references therein) and the literature in $D$-stability is particularly rich, there is still lack of conditions of this kind.
 Indeed, there is a number of sufficient $D$-stability conditions for structured matrices (see, e.g. \cite{JOHN1}) and only a few for unstructured (see, e.g., \cite{OLP}). Also, we have a convenient to use criterion for $n = 3$ (see \cite{CA}) and much less convenient criterion for $n = 4$ (see \cite{KL}). To test a full unstructured $5 \times 5$ matrix is still a problem. Finally, there are algorithms for testing matrix diagonal stability (see \cite{KAB}), which implies $D$-stability. But, as it will be shown by examples, our test works for matrices which are not diagonally stable.

Our new approach exploits the basic linear algebra methods, such as matrix determinant expansion. In practice, we just need to know the principal minors of a matrix, thus matrices of small order ($n \leq 4$) can be easily tested without any software. Note that already in 1970th it was known, that {\it $D$-stability of a matrix depends
entirely on the sequence of its principal minors} (see \cite{JOHN5}). Developing necessary and sufficient criterion of $D$-stability of $n \times n$ matrices, based on this approach, is in progress.

Here, as usual, we use the notation $[n]$ for the set of indices $\{1,2, \ \ldots, \ n\}$. Given a set $\alpha = (i_1, \ \ldots, i_k)$, with $1 \leq i_1 < \ldots < i_k \leq n$, we denote $N(\alpha)$ the cardinality of $\alpha$, i.e. the number of elements in $\alpha$ (obviously, $N(\emptyset) = 0$). For a given matrix ${\mathbf A} = \{a_{ij}\}_{i,j = 1}^n$, we denote $A \begin{pmatrix}\alpha \\ \alpha \end{pmatrix}$ (or $A\begin{pmatrix}i_1 & \ldots & i_k \\ i_1 & \ldots & i_k \end{pmatrix}$) the principal minor of $\mathbf A$, which lies on the intersection of the rows and columns with indices from $\alpha = (i_1, \ \ldots, i_k)$.
The main result of the paper is as follows.

\begin{theorem}\label{cor} Let $\mathbf A$ be a stable matrix such that, for some $k$, $1 \leq k \leq n$, $a_{kk} \neq 0$ and its $(n-1) \times (n-1)$ principal submatrix ${\mathbf A}|_{k}$, obtained by deleting the $k$th row and column, is $D$-stable. Then $\mathbf A$ is $D$-stable if for each two sets of indices $\alpha$, $\beta$, which satisfy the conditions:
 \begin{enumerate}
\item[\rm (a)] $\beta \subseteq \alpha \subseteq[n] \setminus\{k\}$, for $N(\alpha) = 1, \ \ldots, \ n-1$,
\item[\rm (b)] $$N(\beta) = \left\{\begin{array}{cc}  0, \ 2, \ \ldots, \ N(\alpha) & \mbox{if $N(\alpha)$ is even},\\  1, \ 3, \ \ldots, \ N(\alpha) & \mbox{if $N(\alpha)$ is odd} \end{array} \right.$$
\end{enumerate}
     the following inequality holds:
\begin{equation}\label{crit1}\sum_{r = 0}^{N(\alpha\setminus \beta)}\left((-1)^{\chi + r}\sum_{\gamma \subseteq \alpha \setminus \beta, N(\gamma) = r}A\begin{pmatrix}\alpha \setminus \gamma \\ \alpha \setminus \gamma \end{pmatrix}\frac{1}{a_{kk}}A\begin{pmatrix}\beta & \gamma & \hat{k}\\ \beta & \gamma & \hat{k}\\  \end{pmatrix}\right) \geq 0, \end{equation}
where $$\chi = \left\{\begin{array}{cc} n-1 -\frac{N(\beta) + N(\alpha)}{2}, & \mbox{if} \ n-1 - N(\alpha) \ \mbox{ is even}; \\
\lfloor\frac{n-1 - N(\beta)}{2}\rfloor+\lfloor\frac{n-1 - N(\alpha)}{2}\rfloor, & \mbox{if} \ n-1 - N(\alpha) \ \mbox{ is odd}
\end{array} \right.$$ and the notation $\hat{k}$ means that the index $k$ is placed according to its lexicographic order.
\end{theorem}

Note that determinantal inequalities \eqref{crit1} define a class of matrices which is not included in
any known verifiable sufficient condition for $D$-stability.

 The paper is organized as follows. Section 2 collects some basic facts of matrix theory: determinant expansions and Schur complement properties. Section 3 presents the results on $D$-stability: equivalent characteristics, necessary conditions and basic properties. Section 4 is concerned with the proof of the main results.
 In Section 5, the computational aspects are considered and supported by numerical
examples for the dimensions $n = 3$ and $n = 4$. Section 6 contains brief conclusions.

 \section{Preliminary results}
 \subsection{Determinant expansions}
Let us represent an $n \times n$ complex matrix ${\mathbf A} = \{a_{ij}\}_{i,j = 1}^n$, in the following block form
\begin{equation}\label{ym2}
{\mathbf A} = \begin{pmatrix} {\mathbf A}|_{n} & \overline{{\mathbf a}}_{1n} \\
\underline{{\mathbf a}}_{n1} & a_{nn} \end{pmatrix},
\end{equation}
where ${\mathbf A}|_{n}$ is the $(n-1) \times (n-1)$ principal submatrix of $\mathbf A$, obtained by deleting the $n$th row and the $n$th column, $\overline{{\mathbf a}}_{1n} = (a_{1n}, \ \ldots, \ a_{n-1 n})^T$, $\underline{{\mathbf a}}_{n1} = (a_{n1}, \ \ldots, \ a_{n n-1})$.
Later, we shall use the following matrix determinant expansion (see, for example, \cite{BER}, p. 133, Fact 2.14.2).
\begin{lemma}\label{1}
Let an $n \times n$ complex matrix ${\mathbf A} = \{a_{ij}\}_{i,j = 1}^n$ be represented in Form \eqref{ym2}. Then
$$\det{\mathbf A} = \left\{\begin{array}{ccc}\det({\mathbf A}|_{n})(a_{nn} - \underline{{\mathbf a}}_{n1}({\mathbf A}|_{n})^{-1}\overline{{\mathbf a}}_{1n}) & \mbox{if} \ \det({\mathbf A}|_{n}) \neq 0; \\
a_{nn}\det({\mathbf A}|_{n} - \frac{1}{a_{nn}}\overline{{\mathbf a}}_{n1}\underline{{\mathbf a}}_{1n}) & \mbox{if} \ a_{nn} \neq 0.\end{array}\right.$$
\end{lemma}

\subsection{Schur complements}
Given an $n \times n$ complex matrix ${\mathbf A} = \{a_{ij}\}_{i,j = 1}^n$, we decompose it as given in \eqref{ym2} assuming that $a_{nn} \neq 0$. Recall, that the {\it Schur complement} ${\mathbf A}|_{a_{nn}}$ of $a_{nn}$ in $\mathbf A$ is defined as follows:
\begin{equation}\label{sh1}{\mathbf A}|_{a_{nn}}  = {\mathbf A}|_n - \frac{1}{a_{nn}}\overline{{\mathbf a}}_{1n}\underline{{\mathbf a}}_{n1}.\end{equation}

For the entries of the Schur complement, the following equalities hold:
\begin{equation}\label{sh2}{\mathbf A}|_{a_{nn}} = \{b_{ij}\}_{i,j = 1}^{n-1}, \qquad \mbox{where} \qquad b_{ij} = \frac{1}{a_{nn}}A\begin{pmatrix}i & n \\ j & n \end{pmatrix}.\end{equation}

For the determinant of the Schur complement, we have:
\begin{equation}\label{sh3}\det({\mathbf A}|_{a_{nn}}) = \frac{\det({\mathbf A})}{a_{nn}}.\end{equation}

Using the Sylvester determinant identity (see, for example, \cite{BER}, p. 132, Fact 2.14.1), we easily deduce the following formula for the principal minors of the Schur complement ${\mathbf B}:={\mathbf A}|_{a_{nn}}$:
$$B\begin{pmatrix}i_1 &  \ldots & i_k \\ i_1 &  \ldots & i_k\end{pmatrix} = \frac{1}{a_{nn}}A\begin{pmatrix}i_1 &  \ldots & i_k & n \\ i_1 &  \ldots & i_k & n\end{pmatrix},$$
for all $(i_1, \ \ldots, \ i_k)$, $1 \leq i_1 < \ldots < i_k \leq n-1$, and all $k = 1, \ \ldots, \ n-1$.

\section{Basic results on $D$-stablity}
First, recall the following basic results (see, e.g. \cite{JOHN1}).

\begin{lemma}[Elementary properties of $D$-stable matrices]\label{el} If $\mathbf A$ is $D$-stable then $\mathbf A$ is nonsingular and each of the following matrices are also $D$-stable:
\begin{enumerate}
\item[\rm 1.] ${\mathbf A}^T$
\item[\rm 2.] ${\mathbf A}^{-1}$
\item[\rm 3.] ${\mathbf P}^T{\mathbf A}{\mathbf P},$ where $\mathbf P$ is any permutation matrix.
\item[\rm 4.] ${\mathbf D}{\mathbf A}{\mathbf E}$, where $\mathbf D$, $\mathbf E$ are positive diagonal matrices.
\end{enumerate}
\end{lemma}

Here, we recall the following equivalent conditions of matrix $D$-stability (see \cite{JOHN5}, p. 89, Corollary 2).

\begin{theorem}\label{stabcond} Let ${\mathbf A} \in {\mathcal M}^{n \times n}$ be stable. Then the following conditions are equivalent.
\begin{enumerate}
\item[\rm (i)] ${\mathbf A}$ is $D$-stable.
\item[\rm (ii)] $\det({\mathbf A} \pm i{\mathbf D}) \neq 0$ for every positive diagonal matrix $\mathbf D$.
\end{enumerate}
\end{theorem}
The condition (ii) shows that ${\mathbf D}{\mathbf A}$ has no eigenvalues on the imaginary axis for any positive diagonal matrix $\mathbf D$.

Let us recall the following definitions from matrix theory.

{\bf Definition 4.} A matrix ${\mathbf A} \in {\mathcal M}^{n \times n}$ is called:
\begin{enumerate}
\item[-] a {\it $P$-matrix ($P_0$-matrix)} if all its principal minors are positive (respectively, nonnegative), i.e the inequality $A \left(\begin{array}{ccc}i_1 & \ldots & i_k \\ i_1 & \ldots & i_k \end{array}\right) > 0$ (respectively, $\geq 0$)
holds for all $(i_1, \ \ldots, \ i_k), \ 1 \leq i_1 < \ldots < i_k \leq n$, and all $k, \ 1 \leq k \leq n$.
\item[-] a a {\it $P_0^+$-matrix} if it is a $P_0$-matrix and, in addition, the inequality
$$\sum_{(i_1, \ldots, i_k)}A \begin{pmatrix}i_1 & \ldots & i_k \\ i_1 & \ldots & i_k \end{pmatrix} > 0 $$
holds for all $k, \ 1 \leq k \leq n$, i.e. the sums of all principal minors of every fixed order $k$ are positive.
\end{enumerate}

A important necessary for $D$-stability condition was proved in \cite{QR}: {\it if $\mathbf A$ is $D$-stable then $-{\mathbf A}$ is a $P_0^+$-matrix.}

\section{Proof of the main result}

\begin{lemma}\label{explemma} Let an $n\times n$ matrix $\mathbf A$ satisfies the following conditions: $a_{nn} \neq 0$ and $\det({\mathbf A}|_{n} \pm i{\mathbf D}|_{n}) \neq 0$ for every $(n-1)\times(n-1)$ positive diagonal matrix ${\mathbf D}|_{n}$. Then, for each $n \times n$ positive diagonal matrix $\mathbf D$, the following equality holds:
\begin{equation}\label{expan}\det({\mathbf A} \pm i{\mathbf D}) = \pm id_{nn}\det({\mathbf A}|_{n}\pm i{\mathbf D}|_{n}) + a_{nn}\det({\mathbf A}|_{a_{nn}} \pm i {\mathbf D}|_n), \end{equation}
where ${\mathbf D} = \diag\{{\mathbf D}|_n, \ d_{nn}\}$.
\end{lemma}
{\bf Proof.} Let us represent a positive diagonal matrix $\mathbf D$ in the form:
$${\mathbf D} = \diag\{{\mathbf D}|_{n}, d_{nn}\}.$$
Using Expansion \eqref{ym2}, we obtain
$${\mathbf A} \pm i{\mathbf D} = \begin{pmatrix} {\mathbf A}|_{n} & \overline{{\mathbf a}}_{1n} \\
\underline{{\mathbf a}}_{n1} & a_{nn} \end{pmatrix} \pm i\begin{pmatrix} {\mathbf D}|_{n} & \overline{{\mathbf 0}} \\
\underline{{\mathbf 0}} & d_{nn} \end{pmatrix} = \begin{pmatrix} {\mathbf A}|_{n}\pm i{\mathbf D}|_{n} & \overline{{\mathbf a}}_{1n} \\
\underline{{\mathbf a}}_{n1} & a_{nn}\pm id_{nn} \end{pmatrix}.$$

Since $\det({\mathbf A}|_{n} \pm i{\mathbf D}|_{n}) \neq 0$, we use Lemma \ref{1} and expand $\det({\mathbf A} \pm i{\mathbf D})$ as follows:

$$\det({\mathbf A} \pm i{\mathbf D}) = \det({\mathbf A}|_{n}\pm i{\mathbf D}|_{n})(a_{nn} \pm id_{nn} - \underline{{\mathbf a}}_{n1}({\mathbf A}|_{n}\pm i{\mathbf D}|_{n})^{-1}\overline{{\mathbf a}}_{1n})= $$
$$\pm id_{nn}\det({\mathbf A}|_{n}\pm i{\mathbf D}|_{n}) + \det({\mathbf A}|_{n}\pm i{\mathbf D}|_{n})(a_{nn} - \underline{{\mathbf a}}_{n1}({\mathbf A}|_{n}\pm i{\mathbf D}|_{n})^{-1}\overline{{\mathbf a}}_{1n} ) = \ldots$$
Using Lemma \ref{1} again, we get $$\det({\mathbf A}|_{n}\pm i{\mathbf D}|_{n})(a_{nn} - \underline{{\mathbf a}}_{n1}({\mathbf A}|_n + i {\mathbf D}|_{n})^{-1}\overline{{\mathbf a}}_{1n}) = $$
 $$a_{nn}\det({\mathbf A}|_n + i {\mathbf D}|_{n} - \frac{1}{a_{nn}}\underline{{\mathbf a}}_{n1}\overline{{\mathbf a}}_{1n}) = a_{nn}\det({\mathbf A}|_{a_{nn}} \pm i {\mathbf D}|_n).$$
 Thus $$\ldots = \pm id_{nn}\det({\mathbf A}|_{n}\pm i{\mathbf D}|_{n}) + a_{nn}\det({\mathbf A}|_{a_{nn}} \pm i {\mathbf D}|_n). $$ $\square$

 For $k = 1, \ \ldots, \ n$, let us consider the set of indices $(i_1,\ldots,i_{k}),$ $1 \leq i_1 < \ldots < i_k \leq n$. Then denote $\alpha = (j_1, \ \ldots, \ j_l) := [n]\setminus (i_1,\ldots,i_{k})$, here $N(\alpha) = l = n-k$. When $(i_1,\ldots,i_{k}) = [n]$, we obtain $\alpha =\emptyset$ and we set an "empty" principal minor $A\begin{pmatrix}\emptyset\\ \emptyset \end{pmatrix} := 1$.  When $(i_1,\ldots,i_{k}) = \emptyset$, we have $A\begin{pmatrix} \alpha \\ \alpha\end{pmatrix} = A\begin{pmatrix} 1 & \ldots & n \\ 1 & \ldots & n\end{pmatrix}= \det({\mathbf A})$. Let us recall the following lemma from \cite{JOHN5} (see \cite{JOHN5}, p. 89, Lemma 1).
\begin{lemma}\label{john} For a matrix ${\mathbf A} \in {\mathcal M}^{n \times n}$ and an $n \times n$ real diagonal matrix $\mathbf D$, the following formulas hold:
$${\rm Re}(\det({\mathbf A}+ i{\mathbf D})) = \sum_{ p =0}^{ \lfloor\frac{n}{2}\rfloor}\left((-1)^p\sum_{(i_1,\ldots,i_{2p})}A\begin{pmatrix}\alpha \\ \alpha \end{pmatrix}d_{i_1}\ldots d_{i_{2p}}\right), $$
where $\alpha = [n]\setminus(i_1,\ldots,i_{2p})$, the sum is taken with respect to all the subsets $(i_1,\ldots,i_{2p})$, such that $1 \leq i_1 < \ldots < i_{2p} \leq n$.
$${\rm Im}(\det({\mathbf A}+ i{\mathbf D})) = \sum_{p = 1}^{\lceil\frac{n}{2}\rceil}\left((-1)^{p-1}\sum_{(i_1,\ldots,i_{2p-1})}A\begin{pmatrix}\alpha \\ \alpha \end{pmatrix}d_{i_1}\ldots d_{i_{2p-1}}\right),$$
where $\alpha = [n]\setminus(i_1,\ldots,i_{2p-1})$, the sum is taken with respect to all the subsets $(i_1,\ldots,i_{2p-1})$, such that $1 \leq i_1 < \ldots < i_{2p-1} \leq n$.
\end{lemma}

Now let us prove the following crucial statement.

\begin{theorem}\label{ma} Let $\mathbf A$ be a stable matrix such that $a_{nn} \neq 0$ and its $(n-1) \times (n-1)$ principal submatrix ${\mathbf A}|_{n}$ is $D$-stable. Then $\mathbf A$ is $D$-stable, if the equation of $n-1$ variables $$F(d_{11}, \ \ldots, \ d_{n-1n-1}) = 0, $$
with $$F(d_{11}, \ \ldots, \ d_{n-1n-1}) = \left(\sum_{ p =0}^{ \lfloor\frac{n-1}{2}\rfloor}\left((-1)^p\sum_{(i_1,\ldots,i_{2p})}A\begin{pmatrix}\alpha \\ \alpha \end{pmatrix}d_{i_1}\ldots d_{i_{2p}}\right)\right)*$$ $$\left(\sum_{ r =0}^{ \lfloor\frac{n-1}{2}\rfloor}\left((-1)^r\sum_{(k_1,\ldots,k_{2r})}\frac{1}{a_{nn}}A\begin{pmatrix}\beta & n\\ \beta & n \end{pmatrix}d_{k_1}\ldots d_{k_{2r}}\right)\right) + $$
   $$\left(\sum_{l = 1}^{\lceil\frac{n-1}{2}\rceil}\left((-1)^{l-1}\sum_{(j_1,\ldots,j_{2l-1})}A\begin{pmatrix}\gamma \\ \gamma \end{pmatrix}d_{j_1}\ldots d_{j_{2l-1}}\right)\right)* $$ $$\left(\sum_{m=1}^{\lceil\frac{n-1}{2}\rceil}\left((-1)^{m-1}\sum_{(q_1,\ldots,q_{2m-1})}\frac{1}{a_{nn}}A\begin{pmatrix}\omega & n\\ \omega& n \end{pmatrix}d_{q_1}\ldots d_{q_{2m-1}}\right)\right),$$
   where $\alpha=[n-1]\setminus(i_1,\ldots,i_{2p})$, $\beta =[n-1]\setminus(k_1,\ldots,k_{2r})$, $\gamma = [n-1]\setminus(j_1,\ldots,j_{2l-1})$, $\omega = [n-1]\setminus (q_1,\ldots,q_{2m-1})$, has no positive solutions \linebreak $d_{11}, \ \ldots, \ d_{n-1n-1}$.
\end{theorem}
{\bf Proof.}
 Applying Theorem \ref{stabcond} to $\mathbf A$, we obtain that $\mathbf A$ is $D$-stable if and only if $\det({\mathbf A}\pm i{\mathbf D}) \neq 0$. Assume the opposite: let $\det({\mathbf A}+ i{\mathbf D}_0) = 0$ for some positive diagonal matrix ${\mathbf D}_0 = \{d^0_{11}, \ \ldots, \ d^0_{nn}\}$. Since ${\mathbf A}|_{n}$ is $D$-stable, we apply Theorem \ref{stabcond} and obtain that $\det({\mathbf A}|_{n}\pm i{\mathbf D}|_{n}) \neq 0$ for every positive diagonal matrix ${\mathbf D}|_{n}$. Since $\det({\mathbf A}|_{n}+ i{\mathbf D}_0|_{n}) \neq 0$ as well, we obtain by Formula \eqref{expan}: $$\det({\mathbf A} + i{\mathbf D}_0) = id^0_{nn}\det({\mathbf A}|_{n}+ i{\mathbf D}_0|_{n}) + a_{nn}\det({\mathbf A}|_{a_{nn}} + i {\mathbf D}_0|_n). $$ Thus $\det({\mathbf A}+ i{\mathbf D}_0) = 0$ is equivalent to the equation $$id^0_{nn}\det({\mathbf A}|_{n}+ i{\mathbf D}_0|_{n}) + a_{nn}\det({\mathbf A}|_{a_{nn}} + i {\mathbf D}_0|_n) = 0.$$
Taking the equalities between real and imaginary parts, we obtain the following equivalent system:
  $$\left\{\begin{array}{cc} -d^0_{nn}{\rm Im}(\det({\mathbf A}|_{n}+ i{\mathbf D}_0|_{n}))+ a_{nn}{\rm Re}(\det({\mathbf A}|_{a_{nn}}) + i {\mathbf D}_0|_n) = 0; \\
    d^0_{nn}{\rm Re}(\det({\mathbf A}|_{n}+ i{\mathbf D}_0|_{n}))  + a_{nn}{\rm Im}(\det({\mathbf A}|_{a_{nn}}) + i {\mathbf D}_0|_n)=0.\end{array}\right.$$
   The inequality $\det({\mathbf A}|_{n}+ i{\mathbf D}_0|_{n}) \neq 0$ implies that if ${\rm Im}(\det({\mathbf A}|_{n}+ i{\mathbf D}_0|_{n})) = 0$ then ${\rm Re}(\det({\mathbf A}|_{n}+ i{\mathbf D}_0|_{n})) \neq 0$ and conversely. Thus we can extract $d^0_{nn}$ from either the first or the second equation:
     $$d^0_{nn} = \frac{a_{nn}{\rm Re}(\det({\mathbf A}|_{a_{nn}}) + i {\mathbf D}_0|_n)}{{\rm Im}(\det({\mathbf A}|_{n}+ i{\mathbf D}_0|_{n}))} \qquad \mbox{or} \qquad d^0_{nn} = -\frac{a_{nn}{\rm Im}(\det({\mathbf A}|_{a_{nn}}) + i {\mathbf D}_0|_n)}{{\rm Re}(\det({\mathbf A}|_{n}+ i{\mathbf D}_0|_{n}))}. $$
    Substituting $d^0_{nn}$ into the other equation, we get
    $${\rm Re}(\det({\mathbf A}|_{a_{nn}} + i {\mathbf D}_0|_n)){\rm Re}(\det({\mathbf A}|_{n}+ i{\mathbf D}_0|_{n}))  + $$ $$ {\rm Im}(\det({\mathbf A}|_{a_{nn}} + i {\mathbf D}_0|_n)){\rm Im}(\det({\mathbf A}|_{n}+ i{\mathbf D}_0|_{n}))=0 $$

Applying Lemma \ref{john} to both ${\mathbf A}|_{n}$ and ${\mathbf A}|_{a_{nn}}$, we obtain:

  $${\rm Re}(\det({\mathbf A}|_{n}+ i{\mathbf D}_0|_{n})) = \sum_{ p =0}^{ \lfloor\frac{n-1}{2}\rfloor}\left((-1)^p\sum_{(i_1,\ldots,i_{2p})}A\begin{pmatrix}\alpha \\ \alpha \end{pmatrix}d^0_{i_1}\ldots d^0_{i_{2p}}\right),$$
  where $\alpha=[n-1]\setminus(i_1,\ldots,i_{2p})$, the sum is taken with respect to all the subsets $(i_1,\ldots,i_{2p})$, such that $1 \leq i_1 < \ldots < i_{2p} \leq n-1$.
   $${\rm Re}(\det({\mathbf A}|_{a_{nn}}+ i{\mathbf D}_0|_{n})) = \sum_{ r =0}^{ \lfloor\frac{n-1}{2}\rfloor}\left((-1)^r\sum_{(k_1,\ldots,k_{2r})}B\begin{pmatrix}\beta \\ \beta \end{pmatrix}d^0_{k_1}\ldots d^0_{k_{2r}}\right),$$
   where $\beta =[n-1]\setminus(k_1,\ldots,k_{2r})$, the sum is taken with respect to all the subsets $(k_1,\ldots,k_{2r})$, such that $1 \leq k_1 < \ldots < k_{2r} \leq n-1$.
   $${\rm Im}(\det({\mathbf A}|_{n}+ i{\mathbf D}_0|_{n})) = \sum_{l = 1}^{\lceil\frac{n-1}{2}\rceil}\left((-1)^{l-1}\sum_{(j_1,\ldots,j_{2l-1})}A\begin{pmatrix}\gamma \\ \gamma \end{pmatrix}d^0_{j_1}\ldots d^0_{j_{2l-1}}\right),$$
   where $\gamma = [n-1]\setminus(j_1,\ldots,j_{2l-1})$, the sum is taken with respect to all the subsets $(j_1,\ldots,j_{2l-1})$, such that $1 \leq j_1 < \ldots < j_{2l-1} \leq n-1$.
   $${\rm Im}(\det({\mathbf A}|_{a_{nn}}+ i{\mathbf D}_0|_{n})) = \sum_{m=1}^{\lceil\frac{n-1}{2}\rceil}\left((-1)^{m-1}\sum_{(q_1,\ldots,q_{2m-1})}B\begin{pmatrix}\omega \\ \omega \end{pmatrix}d^0_{q_1}\ldots d^0_{q_{2m-1}}\right),$$
   where $\omega = [n-1]\setminus (q_1,\ldots,q_{2m-1})$, the sum is taken with respect to all the subsets $(q_1,\ldots,q_{2m-1})$, such that $1 \leq q_1 < \ldots < q_{2m-1} \leq n-1$.

   Taking into account that
   $$B\begin{pmatrix}i_1 &  \ldots & i_k \\ i_1 &  \ldots & i_k\end{pmatrix} = \frac{1}{a_{nn}}A\begin{pmatrix}i_1 &  \ldots & i_k & n \\ i_1 &  \ldots & i_k & n\end{pmatrix},$$
for all $(i_1, \ \ldots, \ i_k)$, $1 \leq i_1 < \ldots < i_k \leq n-1$, and all $k = 1, \ \ldots, \ n-1$,
  we obtain
   $$\left(\sum_{ p =0}^{ \lfloor\frac{n-1}{2}\rfloor}\left((-1)^p\sum_{(i_1,\ldots,i_{2p})}A\begin{pmatrix}\alpha \\ \alpha \end{pmatrix}d^0_{i_1}\ldots d^0_{i_{2p}}\right)\right)*$$ $$\left(\sum_{r =0}^{ \lfloor\frac{n-1}{2}\rfloor}\left((-1)^r\sum_{(k_1,\ldots,k_{2r})}A\begin{pmatrix}\beta & n \\ \beta & n\end{pmatrix}d^0_{k_1}\ldots d^0_{k_{2r}}\right)\right) + $$
   $$\left(\sum_{l = 1}^{\lceil\frac{n-1}{2}\rceil}\left((-1)^{l-1}\sum_{(j_1,\ldots,j_{2l-1})}A\begin{pmatrix}\gamma \\ \gamma \end{pmatrix}d^0_{j_1}\ldots d^0_{j_{2l-1}}\right)\right)* $$ $$\left(\sum_{m=1}^{\lceil\frac{n-1}{2}\rceil}\left((-1)^{m-1}\sum_{(q_1,\ldots,q_{2m-1})}A\begin{pmatrix}\omega & n \\ \omega & n \end{pmatrix}d^0_{q_1}\ldots d^0_{q_{2m-1}}\right)\right)= 0.$$
   The above equality means exactly that the equation $F(d_{11}, \ \ldots, \ d_{n-1n-1}) =0$ has a positive solution $d^0_{11}, \ \ldots, \ d^0_{n-1n-1}$. We obtain a contradiction.
 $\square$

{\bf Proof of Theorem 1.} Without loss the generality of the reasoning, we assume that $k = n$. Otherwise, we shall consider ${\mathbf B} := {\mathbf P} {\mathbf A}{\mathbf P}^T$, where ${\mathbf P}$ is a matrix of the permutation, which interchanges $n$ and $k$. According to Lemma \ref{el}, the matrix ${\mathbf B}$ is $D$-stable if and only if $\mathbf A$ is $D$-stable.

 Consider the function $F(d_{11}, \ \ldots, \ d_{n-1n-1})$. After opening the brackets, it becomes the sum of monomials in variables $d_{11}, \ \ldots, \ d_{n-1n-1}$ of the highest degree $2$. Given two sets $\alpha$, $\beta$ such that $\beta \subseteq \alpha \subseteq[n-1]$, let us show that the expression
 $$\sum_{r = 0}^{N(\alpha\setminus \beta)}\left((-1)^{\chi + r}\sum_{\gamma \subseteq \alpha \setminus \beta, N(\gamma) = r}A\begin{pmatrix}\alpha \setminus \gamma \\ \alpha \setminus \gamma \end{pmatrix}\frac{1}{a_{nn}}A\begin{pmatrix}\beta & \gamma & n\\ \beta & \gamma & n\\  \end{pmatrix}\right)$$ equals the coefficient of some monomial in $F(d_{11}, \ \ldots, \ d_{n-1n-1})$. Indeed, let us define $(i_1, \ \ldots, \ i_k) = [n-1] \setminus \beta$ and $(j_1, \ \ldots, \ j_l) = [n-1] \setminus \alpha$. Obviously, $(j_1, \ \ldots, \ j_l)\subseteq (i_1, \ \ldots, \ i_k) \subseteq [n-1]$. Now let us consider the monomial $d_{i_1}^{t_1}\ldots d_{i_k}^{t_k}$, where $t_r = 2$ if $i_r \in (j_1, \ \ldots, \ j_l)$ and $t_r = 1$ otherwise. To calculate the coefficient of $d_{i_1}^{t_1}\ldots d_{i_k}^{t_k}$, let us count all the decompositions of the form
 $$d_{i_1}^{p_1}\ldots d_{i_k}^{p_k} = (d_{q_1}\ldots d_{q_s})(d_{p_1}\ldots d_{p_m}),$$
 where $s$ and $m$ are of the same oddity. For this, we take all the possible sets $$\gamma \subseteq (i_1, \ \ldots, \ i_k) \setminus (j_1, \ \ldots, \ j_l) = \alpha \setminus \beta,$$
starting from $\gamma = \emptyset$ till $\gamma = \alpha \setminus \beta$. Then, we have
$$(q_1, \ \ldots, \ q_s) = (i_1, \ \ldots, \ i_k) \setminus \gamma = [n-1] \setminus (\beta\cup\gamma);$$
$$(p_1, \ \ldots, \ p_m) = (j_1, \ \ldots, \ j_l) \cup \gamma = [n-1] \setminus (\alpha \setminus\gamma).$$
 Each product $(d_{q_1}\ldots d_{q_s})(d_{p_1}\ldots d_{p_m})$ has the coefficient of the form $$A\begin{pmatrix}\alpha \setminus \gamma \\ \alpha \setminus \gamma \end{pmatrix}\frac{1}{a_{nn}}A\begin{pmatrix}\beta & \gamma & n\\ \beta & \gamma & n\\  \end{pmatrix},$$
 taken with the corresponding sign.
Let us define the sign of the first summand $$A\begin{pmatrix}\alpha  \\ \alpha \end{pmatrix}\frac{1}{a_{nn}}A\begin{pmatrix}\beta  & n\\ \beta  & n\\  \end{pmatrix},$$ which corresponds to $\gamma = \emptyset$.
 If $k = n-1 - N(\beta)$ is even then $l = n-1 - N(\alpha)$ is also even. Then, according to the formula for $F(d_{11}, \ \ldots, \ d_{n-1n-1})$, the minor $A\begin{pmatrix}\alpha  \\ \alpha \end{pmatrix}$ comes with the sign $(-1)^{\frac{k}{2}}$ and $ \frac{1}{a_{nn}}A\begin{pmatrix}\beta  & n\\ \beta  & n\\  \end{pmatrix}$ comes with the sign $(-1)^{\frac{l}{2}}$. Thus their product $A\begin{pmatrix}\alpha  \\ \alpha \end{pmatrix}\frac{1}{a_{nn}}A\begin{pmatrix}\beta  & n\\ \beta  & n\\  \end{pmatrix}$ is involved with the sign $$(-1)^{\frac{k+l}{2}} = (-1)^{n-1 -\frac{N(\beta) + N(\alpha)}{2}}.$$
 In the case when $k$ is odd, we have that $l$ is also odd, and, according to the expression for $F(d_{11}, \ \ldots, \ d_{n-1n-1})$, the product is involved with the sign $$(-1)^{\lfloor\frac{k}{2}\rfloor+\lfloor\frac{l}{2}\rfloor} = (-1)^{\lfloor\frac{n-1 - N(\beta)}{2}\rfloor+\lfloor\frac{n-1 - N(\alpha)}{2}\rfloor}.$$ Now let us show, that if $N(\gamma) = r$ is even, the summand $$A\begin{pmatrix}\alpha \setminus \gamma \\ \alpha \setminus \gamma \end{pmatrix}\frac{1}{a_{nn}}A\begin{pmatrix}\beta & \gamma & n\\ \beta & \gamma & n\\  \end{pmatrix}$$ comes with the same sign. Indeed, in this case, both $s = n-1 - N(\beta) - r$ and $m = n-1 - N(\alpha) + r$ have the same oddity as $k = n-1 - N(\beta)$ and $l = n-1 - N(\alpha)$. In the case when they are even, we obtain the sign $$(-1)^{\frac{s+m}{2}} = (-1)^{n-1 -\frac{N(\beta) - r + r + N(\alpha)}{2}} = (-1)^{n-1 -\frac{N(\beta) + N(\alpha)}{2}}.$$ In the case of odd, $$(-1)^{\lfloor\frac{s}{2}\rfloor+\lfloor\frac{m}{2}\rfloor}=(-1)^{\lfloor\frac{n-1 - N(\beta)-r}{2}\rfloor+\lfloor\frac{n-1 - N(\alpha)+r}{2}\rfloor} = (-1)^{\lfloor\frac{n-1 - N(\beta)}{2}\rfloor+\lfloor\frac{n-1 - N(\alpha)}{2}\rfloor - \frac{r}{2}+ \frac{r}{2}} = $$ $$ (-1)^{\lfloor\frac{n-1 - N(\beta)}{2}\rfloor+\lfloor\frac{n-1 - N(\alpha)}{2}\rfloor},$$ since $r$ is even. In the case, when $r$ is odd, it changes the oddity, and it is easy to see, that the sign changes to the opposite.
 
  Thus, all the coefficients of monomials involved in the function \linebreak $F(d_{11}, \ \ldots, \ d_{n-1n-1})$ are nonnegative.
  By noticing that at least one monomial is involved with the positive coefficient, we complete the proof. $\square$
\section{Examples}

\subsection{Case of matrices $3 \times 3$}
 Consider ${\mathbf A} = \{a_{ij}\}_{i,j = 1}^3$. Let us check its $D$-stability, using Theorem \ref{cor}. Here, we assume $k = 3$, that means, $a_{33} \neq 0$ and the principal submatrix ${\mathbf A}|_3 = \{a_{ij}\}_{i,j = 1}^2$ is $D$-stable. First, let us check stability of $\mathbf A$. The well-known Routh--Hurwitz criterion gives us $$\det({\mathbf A}) < 0; \qquad a_{11} + a_{22} + a_{33} < 0; \qquad A\begin{pmatrix} 1 & 2 \\ 1 & 2 \end{pmatrix} + A\begin{pmatrix} 1 & 3 \\ 1 & 3 \end{pmatrix} + A\begin{pmatrix} 2 & 3 \\ 2 & 3 \end{pmatrix} > 0;$$ and
$$(a_{11} + a_{22} + a_{33})\left(A\begin{pmatrix} 1 & 2 \\ 1 & 2 \end{pmatrix} + A\begin{pmatrix} 1 & 3 \\ 1 & 3 \end{pmatrix} + A\begin{pmatrix} 2 & 3 \\ 2 & 3 \end{pmatrix} \right) < \det({\mathbf A}).$$

Assuming $D$-stability of submatrices, we additionally get $a_{11} \leq 0$, $a_{22} \leq 0$, $a_{11} + a_{22} < 0$, $A\begin{pmatrix} 1 & 2 \\ 1 & 2 \end{pmatrix} > 0$.

Conditions \eqref{crit1} give four more inequalities:

$$\alpha = \{1\}, \beta=\{1\}: \qquad a_{11}\frac{1}{a_{33}}A\begin{pmatrix} 1 & 3 \\ 1 & 3 \end{pmatrix} \geq 0;$$
$$\alpha = \{2\}, \beta=\{2\}: \qquad a_{22}\frac{1}{a_{33}}A\begin{pmatrix} 2 & 3 \\ 2 & 3 \end{pmatrix} \geq 0; $$
$$\alpha = \{1,2\}, \beta=\{\emptyset\}: $$ $$  -A\begin{pmatrix} 1 & 2 \\ 1 & 2 \end{pmatrix} + \frac{a_{11}}{a_{33}}A\begin{pmatrix} 2 & 3 \\ 2 & 3 \end{pmatrix} + \frac{a_{22}}{a_{33}}A\begin{pmatrix} 1 & 3 \\ 1 & 3 \end{pmatrix} - \frac{1}{a_{33}}\det({\mathbf A}) \geq 0;$$
$$\alpha = \{1,2\}, \beta=\{1,2\}: \qquad A\begin{pmatrix} 1 & 2 \\ 1 & 2 \end{pmatrix}\frac{1}{a_{33}}\det({\mathbf A}) \geq 0. $$

Applying the above reasoning for $k = 1$ and $k = 2$, we obtain other conditions, sufficient for $D$-stability. Summarizing, we obtain the following statement.
\begin{proposition}
Given a stable $3 \times 3$ matrix $\mathbf A$, such that $-{\mathbf A}$ is a $P$-matrix. Then $\mathbf A$ is $D$-stable if at least one of the following inequalities holds:
\begin{equation}\label{condn3}-A\begin{pmatrix} 1 & 2 \\ 1 & 2 \end{pmatrix} + \frac{a_{11}}{a_{33}}A\begin{pmatrix} 2 & 3 \\ 2 & 3 \end{pmatrix} + \frac{a_{22}}{a_{33}}A\begin{pmatrix} 1 & 3 \\ 1 & 3 \end{pmatrix} - \frac{1}{a_{33}}\det({\mathbf A}) \geq 0; \end{equation}
\begin{equation}\label{condn4}-A\begin{pmatrix} 1 & 3 \\ 1 & 3 \end{pmatrix} + \frac{a_{11}}{a_{22}}A\begin{pmatrix} 2 & 3 \\ 2 & 3 \end{pmatrix} + \frac{a_{33}}{a_{22}}A\begin{pmatrix} 1 & 2 \\ 1 & 2 \end{pmatrix} - \frac{1}{a_{22}}\det({\mathbf A}) \geq 0; \end{equation}
\begin{equation}\label{condn5}-A\begin{pmatrix} 2 & 3 \\ 2 & 3 \end{pmatrix} + \frac{a_{22}}{a_{11}}A\begin{pmatrix} 1 & 3 \\ 1 & 3 \end{pmatrix} + \frac{a_{33}}{a_{11}}A\begin{pmatrix} 1 & 2 \\ 1 & 2 \end{pmatrix} - \frac{1}{a_{11}}\det({\mathbf A}) \geq 0. \end{equation}
\end{proposition}

The proof follows from the above reasoning.

{\bf Example 1.} Let us show, that Conditions \eqref{crit1} are not necessary for $D$-stability. For this, let us consider the following example of $3 \times 3$ $D$-stable matrix from \cite{JOHN4}:
$${\mathbf A} = \begin{pmatrix} -6 & -5 & 1 \\ -1 & -2 & -5 \\ -5 & 3 & -1 \end{pmatrix}.$$
In this case, $$A\begin{pmatrix} 1 & 2 \\ 1 & 2 \end{pmatrix} = 7; \qquad A\begin{pmatrix} 1 & 3 \\ 1 & 3 \end{pmatrix} = 11; \qquad A\begin{pmatrix} 2 & 3 \\ 2 & 3 \end{pmatrix} = 17; \qquad \det({\mathbf A}) = - 235.$$
However, Condition \eqref{condn3} gives us:
$$17*6 + 11*2 - 235 - 7 = - 118 < 0. $$
Checking Conditions \eqref{condn4} and \eqref{condn5}, we also obtain negative values.

{\bf Example 2.} Let $${\mathbf A} = \begin{pmatrix} -1 & 0 & q \\
-1 & -1 & 0 \\ -1 & -1 & -1 \end{pmatrix}.$$ This matrix is known to be stable for $q > -\frac{8}{3}$. Also note, that for $q = 50$ it is shown to be {\bf not} diagonally stable (see \cite{HART}, p. 202).
In this case, $$A\begin{pmatrix} 1 & 2 \\ 1 & 2 \end{pmatrix} = A\begin{pmatrix} 2 & 3 \\ 2 & 3 \end{pmatrix} = 1; \qquad A\begin{pmatrix} 1 & 3 \\ 1 & 3 \end{pmatrix} = 1 + q;  \qquad \det({\mathbf A}) = - 1.$$
By Condition \eqref{condn3}, we obtain the following range for parameter $q$, which would guarantee $D$-stability:
$$q \geq -1 \qquad \mbox{and} \qquad -1 + (1 + q) + 1 -1 \geq 0 \ \Rightarrow \ q \geq 0.$$
Condition \eqref{condn4} gives an additional interval:
$$q \geq -1 \qquad \mbox{and} \qquad -(1 + q) + 1 +1 - 1 \geq 0 \ \Rightarrow \ q \in [-1; 0].$$
 Thus, sufficient for $D$-stability Conditions \eqref{crit1} establishes $D$-stability of $\mathbf A$ for all $q \geq -1$ independently from its diagonal stability, that coincides with the results from \cite{IJP} and \cite{JOHNT}, obtained by different methods.
\subsection{Case of matrices $4 \times 4$} Now let ${\mathbf A} = \{a_{ij}\}_{i,j = 1}^4$. The Routh--Hurwitz criterion gives us $${\rm Tr}(\mathbf A) < 0; \quad \sum_{1 \leq i < j \leq 4}A\begin{pmatrix} i & j \\ i & j \end{pmatrix}  > 0; \quad \sum_{1 \leq i < j < k \leq 4}A\begin{pmatrix} i & j & k \\ i & j & k \end{pmatrix}  < 0; \quad \det({\mathbf A}) > 0$$ and
$$\left({\rm Tr}(\mathbf A)\right)\left(\sum_{1 \leq i < j \leq 4}A\begin{pmatrix} i & j \\ i & j \end{pmatrix} \right)\left(\sum_{1 \leq i < j < k \leq 4}A\begin{pmatrix} i & j & k \\ i & j & k \end{pmatrix} \right) > $$ $$\left(\sum_{1 \leq i < j < k \leq 4}A\begin{pmatrix} i & j & k \\ i & j & k \end{pmatrix} \right)^2 + \left( {\rm Tr}(\mathbf A)\right)^2\det({\mathbf A}).$$

Assuming $D$-stability of the $3 \times 3$ leading principal submatrix ${\mathbf A}|_4$, we obtain that $-{\mathbf A}|_4$ is a $P_0^+$-matrix. Note, that $D$-stability of a submatrix may be established by any possible way, it may not satisfy Conditions \eqref{crit1} for $n = 3$.

Conditions \eqref{crit1} for $n = 4$, $k = 4$ give 13 more inequalities:

$$\alpha = \{1\}, \beta=\{1\}: \qquad a_{11}\frac{1}{a_{44}}A\begin{pmatrix} 1 & 4 \\ 1 & 4 \end{pmatrix} \geq 0;$$
$$\alpha = \{2\}, \beta=\{2\}: \qquad a_{22}\frac{1}{a_{44}}A\begin{pmatrix} 2 & 4 \\ 2 & 4 \end{pmatrix} \geq 0; $$
$$\alpha = \{3\}, \beta=\{3\}: \qquad a_{33}\frac{1}{a_{44}}A\begin{pmatrix} 3 & 4 \\ 3 & 4 \end{pmatrix} \geq 0; $$
$$\alpha = \{1,2\}, \beta=\{\emptyset\}: $$ $$  -A\begin{pmatrix} 1 & 2 \\ 1 & 2 \end{pmatrix} + \frac{a_{11}}{a_{44}}A\begin{pmatrix} 2 & 4 \\ 2 & 4 \end{pmatrix} + \frac{a_{22}}{a_{44}}A\begin{pmatrix} 1 & 4 \\ 1 & 4 \end{pmatrix} - \frac{1}{a_{44}}A\begin{pmatrix} 1 & 2 & 4 \\ 1 & 2 & 4 \end{pmatrix} \geq 0;$$
$$\alpha = \{1,2\}, \beta=\{1,2\}: \qquad A\begin{pmatrix} 1 & 2 \\ 1 & 2 \end{pmatrix}\frac{1}{a_{44}}A\begin{pmatrix} 1 & 2 & 4 \\ 1 & 2 & 4 \end{pmatrix} \geq 0; $$
$$\alpha = \{1,3\}, \beta=\{\emptyset\}: $$ $$  -A\begin{pmatrix} 1 & 3 \\ 1 & 3 \end{pmatrix} + \frac{a_{11}}{a_{44}}A\begin{pmatrix} 3 & 4 \\ 3 & 4 \end{pmatrix} + \frac{a_{33}}{a_{44}}A\begin{pmatrix} 1 & 4 \\ 1 & 4 \end{pmatrix} - \frac{1}{a_{44}}A\begin{pmatrix} 1 & 3 & 4 \\ 1 & 3 & 4 \end{pmatrix} \geq 0;$$
$$\alpha = \{1,3\}, \beta=\{1,3\}: \qquad A\begin{pmatrix} 1 & 3 \\ 1 & 3 \end{pmatrix}\frac{1}{a_{44}}A\begin{pmatrix} 1 & 3 & 4 \\ 1 & 3 & 4 \end{pmatrix} \geq 0; $$
$$\alpha = \{2,3\}, \beta=\{\emptyset\}: $$ $$  -A\begin{pmatrix} 2 & 3 \\ 2 & 3 \end{pmatrix} + \frac{a_{33}}{a_{44}}A\begin{pmatrix} 2 & 4 \\ 2 & 4 \end{pmatrix} + \frac{a_{22}}{a_{44}}A\begin{pmatrix} 3 & 4 \\ 3 & 4 \end{pmatrix} - \frac{1}{a_{44}}A\begin{pmatrix} 2 & 3 & 4 \\ 2 & 3 & 4 \end{pmatrix} \geq 0;$$
$$\alpha = \{2,3\}, \beta=\{2,3\}: \qquad A\begin{pmatrix} 2 & 3 \\ 2 & 3 \end{pmatrix}\frac{1}{a_{44}}A\begin{pmatrix} 2 & 3 & 4 \\ 2 & 3 & 4 \end{pmatrix} \geq 0; $$
Then, multiplying by $a_{44}$ and taking into account that $a_{44} < 0$, we get:
$$\alpha = \{1,2,3\}, \beta=\{1\}: $$ $$ -A\begin{pmatrix} 1 & 2 & 3 \\ 1 & 2 & 3 \end{pmatrix}A\begin{pmatrix} 1 & 4 \\ 1 & 4 \end{pmatrix} + A\begin{pmatrix} 1 & 3 \\ 1 & 3 \end{pmatrix}A\begin{pmatrix} 1 & 2 & 4 \\ 1 & 2 & 4 \end{pmatrix} + $$ $$ A\begin{pmatrix} 1 & 2 \\ 1 & 2 \end{pmatrix}A\begin{pmatrix} 1 & 3 & 4 \\ 1 & 3 & 4 \end{pmatrix} - a_{11}\det({\mathbf A}) \leq 0; $$
$$\alpha = \{1,2,3\}, \beta=\{2\}: $$ $$ -A\begin{pmatrix} 1 & 2 & 3 \\ 1 & 2 & 3 \end{pmatrix}A\begin{pmatrix} 2 & 4 \\ 2 & 4 \end{pmatrix} + A\begin{pmatrix} 2 & 3 \\ 2 & 3 \end{pmatrix}A\begin{pmatrix} 1 & 2 & 4 \\ 1 & 2 & 4 \end{pmatrix} +  $$ $$ A\begin{pmatrix} 1 & 2 \\ 1 & 2 \end{pmatrix}A\begin{pmatrix} 2 & 3 & 4 \\ 2 & 3 & 4 \end{pmatrix} - a_{22}\det({\mathbf A}) \leq 0; $$
$$\alpha = \{1,2,3\}, \beta=\{3\}: $$ $$ -A\begin{pmatrix} 1 & 2 & 3 \\ 1 & 2 & 3 \end{pmatrix}A\begin{pmatrix} 3 & 4 \\ 3 & 4 \end{pmatrix} + A\begin{pmatrix} 1 & 3 \\ 1 & 3 \end{pmatrix}A\begin{pmatrix} 2 & 3 & 4 \\ 2 & 3 & 4 \end{pmatrix} + $$ $$ A\begin{pmatrix} 2 & 3 \\ 2 & 3 \end{pmatrix}A\begin{pmatrix} 1 & 3 & 4 \\ 1 & 3 & 4 \end{pmatrix} - a_{33}\det({\mathbf A}) \leq 0; $$

$$\alpha = \{1,2,3\}, \beta=\{1,2,3\}: \qquad A\begin{pmatrix} 1 & 2 & 3 \\ 1 & 2 & 3 \end{pmatrix}\frac{1}{a_{44}}\det{\mathbf A} \geq 0. $$

Thus, for a stable $4 \times 4$ $P$-matrix with to ensure $D$-stability using Conditions \eqref{crit1} with $k = 4$, it is enough to check six inequalities: three inequalities of the following form:
\begin{equation}\label{two} \frac{a_{ii}}{a_{44}}A\begin{pmatrix} j & 4 \\ j & 4 \end{pmatrix} + \frac{a_{jj}}{a_{44}}A\begin{pmatrix} i & 4 \\ i & 4 \end{pmatrix} \geq A\begin{pmatrix} i & j \\ i & j \end{pmatrix} + \frac{1}{a_{44}}A\begin{pmatrix} i & j & 4 \\ i & j & 4 \end{pmatrix}, \qquad 1 \leq i < j \leq 3; \end{equation}
three inequalities of the following form
\begin{equation}\label{one}A\begin{pmatrix} j & i \\ j & i \end{pmatrix}A\begin{pmatrix} k & i & 4 \\ k & i & 4 \end{pmatrix} + A\begin{pmatrix} k & i \\ k & i \end{pmatrix}A\begin{pmatrix} j & i & 4 \\ j & i & 4 \end{pmatrix} \leq $$ $$ A\begin{pmatrix} 1 & 2 & 3 \\ 1 & 2 & 3 \end{pmatrix}A\begin{pmatrix} i & 4 \\ i & 4 \end{pmatrix} + a_{ii}\det({\mathbf A}),\end{equation}
where $i = 1,2,3$, $(j,k) = (1,2,3) \setminus\{i\}$.

Using Conditions \eqref{crit1} with other possible values of $k$ $(k = 1, \ 2, \ 3)$, we obtain other sets of inequalities that  are sufficient for matrix $D$-stability.

{\bf Example 3.} Consider the following parameter-dependent matrix
 $${\mathbf A} = \begin{pmatrix} -1 & 0 & q & p \\
-1 & -1 & 0 & 0 \\ -1 & -1 & -1 & 0 \\ -1 & -1 & -1 & -1 \end{pmatrix}.$$

As it was shown by Example 2, its $3 \times 3$ leading principal submatrix is $D$-stable for the values of $q \geq -1$.

Calculating the principal minors of $\mathbf A$, we get:
$$a_{ii} = -1, \qquad i = 1, \ \ldots, \ 4, $$ $$ A\begin{pmatrix} 1 & 2 \\ 1 & 2 \end{pmatrix} = A\begin{pmatrix} 2 & 3 \\ 2 & 3 \end{pmatrix} = A\begin{pmatrix} 2 & 4 \\ 2 & 4 \end{pmatrix} = A\begin{pmatrix} 3 & 4 \\ 3 & 4 \end{pmatrix} = 1;$$
$$A\begin{pmatrix} 1 & 3 \\ 1 & 3 \end{pmatrix} = 1 + q; \qquad A\begin{pmatrix} 1 & 4 \\ 1 & 4 \end{pmatrix} = 1 + p; \qquad \det({\mathbf A}) = 1;$$
$$ A\begin{pmatrix} 1 & 2 & 3 \\ 1 & 2 & 3 \end{pmatrix} = A\begin{pmatrix} 2 & 3 & 4 \\ 2 & 3 & 4 \end{pmatrix} = A\begin{pmatrix} 1 & 2 & 4 \\ 1 & 2 & 4 \end{pmatrix} = - 1; \quad A\begin{pmatrix} 1 & 3 & 4 \\ 1 & 3 & 4 \end{pmatrix} = - q - 1;$$
Establishing stability range by the Routh--Hurwitz criterion, we obtain:
 $${\rm Tr}(\mathbf A) = -4 < 0; \quad \sum_{1 \leq i < j \leq 4}A\begin{pmatrix} i & j \\ i & j \end{pmatrix} = 6 + (p+q)  > 0; $$ $$ \sum_{1 \leq i < j < k \leq 4}A\begin{pmatrix} i & j & k \\ i & j & k \end{pmatrix} = -4 - q  < 0; \quad \det({\mathbf A}) = 1 > 0$$ and
$$-4(6+(p+q))(-4 -q) > (-4 -q)^2 + 16 \qquad \Rightarrow \qquad p > - \frac{(q+8)(3q + 8)}{4(q+4)}.$$
Summarizing, we obtain that $\mathbf A$ is stable for $q > -4$, $p > -2$, $p > - \frac{(q+8)(3q + 8)}{4(q+4)}$.

Inequalities \eqref{two} give the following conditions:
$$A\begin{pmatrix} 2 & 4 \\ 2 & 4 \end{pmatrix} + A\begin{pmatrix} 1 & 4 \\ 1 & 4 \end{pmatrix} - A\begin{pmatrix} 1 & 2 \\ 1 & 2 \end{pmatrix} + A\begin{pmatrix} 1 & 2 & 4 \\ 1 & 2 & 4 \end{pmatrix} = 1+(1+p)-1-1 = p \geq 0;$$
$$A\begin{pmatrix} 3 & 4 \\ 3 & 4 \end{pmatrix} + A\begin{pmatrix} 1 & 4 \\ 1 & 4 \end{pmatrix} - A\begin{pmatrix} 1 & 3 \\ 1 & 3 \end{pmatrix} + A\begin{pmatrix} 1 & 3 & 4 \\ 1 & 3 & 4 \end{pmatrix} = $$ $$ 1 + (1+p) - (1+q) + (-q-1) = p-2q \geq 0;$$
$$A\begin{pmatrix} 3 & 4 \\ 3 & 4 \end{pmatrix} + A\begin{pmatrix} 2 & 4 \\ 2 & 4 \end{pmatrix} - A\begin{pmatrix} 2 & 3 \\ 2 & 3 \end{pmatrix} + A\begin{pmatrix} 2 & 3 & 4 \\ 2 & 3 & 4 \end{pmatrix} = 1 + 1 -1 -1 = 0;$$

Inequalities \eqref{one} give:

$$ -A\begin{pmatrix} 1 & 2 & 3 \\ 1 & 2 & 3 \end{pmatrix}A\begin{pmatrix} 1 & 4 \\ 1 & 4 \end{pmatrix} + A\begin{pmatrix} 1 & 3 \\ 1 & 3 \end{pmatrix}A\begin{pmatrix} 1 & 2 & 4 \\ 1 & 2 & 4 \end{pmatrix} + $$ $$ A\begin{pmatrix} 1 & 2 \\ 1 & 2 \end{pmatrix}A\begin{pmatrix} 1 & 3 & 4 \\ 1 & 3 & 4 \end{pmatrix} - a_{11}\det({\mathbf A}) = (1+p) - (1+q) +(-q-1) +1 = p - 2q \leq 0; $$
$$ -A\begin{pmatrix} 1 & 2 & 3 \\ 1 & 2 & 3 \end{pmatrix}A\begin{pmatrix} 2 & 4 \\ 2 & 4 \end{pmatrix} + A\begin{pmatrix} 2 & 3 \\ 2 & 3 \end{pmatrix}A\begin{pmatrix} 1 & 2 & 4 \\ 1 & 2 & 4 \end{pmatrix} + $$ $$ A\begin{pmatrix} 1 & 2 \\ 1 & 2 \end{pmatrix}A\begin{pmatrix} 2 & 3 & 4 \\ 2 & 3 & 4 \end{pmatrix} - a_{22}\det({\mathbf A}) =  1 -1 -1 +1 = 0; $$
$$ -A\begin{pmatrix} 1 & 2 & 3 \\ 1 & 2 & 3 \end{pmatrix}A\begin{pmatrix} 3 & 4 \\ 3 & 4 \end{pmatrix} + A\begin{pmatrix} 1 & 3 \\ 1 & 3 \end{pmatrix}A\begin{pmatrix} 2 & 3 & 4 \\ 2 & 3 & 4 \end{pmatrix} + $$ $$ A\begin{pmatrix} 2 & 3 \\ 2 & 3 \end{pmatrix}A\begin{pmatrix} 1 & 3 & 4 \\ 1 & 3 & 4 \end{pmatrix} - a_{33}\det({\mathbf A}) = 1 - (1+q) + (-q-1) +1 = -2q \leq 0. $$
Summarizing the above conditions, we establish $D$-stability of $\mathbf A$ for the values of parameters $p = 2q \geq 0$. Thus, a $4 \times 4$ matrix of the form $${\mathbf A} = \begin{pmatrix} -1 & 0 & q & 2q \\
-1 & -1 & 0 & 0 \\ -1 & -1 & -1 & 0 \\ -1 & -1 & -1 & -1 \end{pmatrix}$$ is $D$-stable for any positive values of $q$.
Applying Conditions \eqref{crit1} with $k = 1, \ 2, \ 3$, we obtain other possible values for the parameters $p$ and $q$.
\section{Conclusions} Applying the determinant expansion formulas, as in the proof of Theorem \ref{ma}, allows us to obtain less and less conservative sufficient for $D$-stability conditions. Getting to the end, we shall obtain the necessary and sufficient criterion of $D$-stability in terms of inequalities between the principal minors of an $n \times n$ matrix. Inequalities \eqref{crit1} are also of interest from the point of view of structured matrices.
\section*{Acknowledgements} The research was supported by the National Science Foundation of China grant number
12050410229.
{}

\begin{thebibliography}{}
\bibitem{AM}
K.J. Arrow, M. McManus, {\it A note on dynamical stability,} Econometrica {\bf 26} (1958), 448-454.

\bibitem{BELL}
R. Bellman, {\it Introduction to matrix analysis}, McGraw Hill, New York, 2nd edition (1970).

\bibitem{BER}
D.S. Bernstein, {\it Matrix mathematics: theory, facts and formulas}, Princeton University Press (2009).

\bibitem{CA}
B. Cain, {\it Real, $3 \times 3$, $D$-stable matrices,} J. Res. Nat. Bur. Standards Sect. B, {\bf 80B} (1976), 75–77.

\bibitem{CROSS}
G.W. Cross, {\it Three types of matrix stability}, Linear Algebra Appl., {\bf 20} (1978), pp. 253--263.

\bibitem{HART}
 D.J. Hartfiel, {\it Concerning the interior of the $D$-stable matrices,} Linear Algebra Appl., {\bf 30} (1980), pp. 201--207.

\bibitem{IJP}
S.T. Impram, R. Johnson, R. Pavani, {\it The $D$-stability problem for $4 \times 4$ real matrices,} Archivum Mathematicum (BRNO) {\bf 41} (2005), 439--450.

\bibitem{JOHN1}
C.R. Johnson, {\it Sufficient conditions for $D$-stability,} Journal of Economic Theory {\bf 9} (1974), 53-62.

\bibitem{JOHN4}
C.R. Johnson, {\it Second, third and fourth order $D$-stability,} J. Research Nat. Bureau Standards USA {\bf B78(1)} (1974), 11-13.

\bibitem{JOHN5}
C.R. Johnson, {\it A characterization of the nonlinearity of $D$-stability,} Journal of Mathematical Economics {\bf 2}, 87-91 (1975).

\bibitem{JOHNT}
R. Johnson, A. Tesi, {\it On the $D$-stability problem for real matrices,} Bollettino dell'Unione Matematica Italiana, {\bf 2-B}, 299--314 (1999).

\bibitem{KL}
G.V. Kanovei and D.O. Logofet, {\it $D$-stability of 4-by-4 matrices}, Comput. Math. Math. Phys., {\bf 38} (1998), pp. 1369--1374.

\bibitem{KAB}
E. Kaszkurewicz, A. Bhaya, {\it Matrix diagonal stability in systems and computation}, Springer (2000).

\bibitem{KU2}
O. Kushel, {\it Unifying matrix stabiity concepts with a view to applications}, SIAM Rev.,{\bf 61(4)} (2019), 643-729

\bibitem{OLP}
R. Oliveira, P. Peres {\it A simple and less conservative test for $D$-stability,} SIAM J. Matrix Anal. Appl., {\bf 26} (2005), 415--425.

\bibitem{QR}
J.P. Quirk, R. Ruppert, {\it Qualitative economics and the stability of equilibrium,} Rev. Econom. Studies {\bf 32} (1965), 311-326.


\end{thebibliography}
\end{document}